\documentclass[10pt,reqno]{amsart}

\usepackage{amsfonts,latexsym}

\title[Boundary singularities for weak solutions]{Boundary singularities for weak solutions of
semilinear elliptic problems}

\author{Manuel del Pino, Monica Musso, Frank Pacard}

\newcommand{\ee}{\varphi_1}

\newcommand{\foral}{\hbox{ for all }}

\newcommand{\RR}{\mathbb{R}}

\newcommand{\ve}{\varepsilon}
\newcommand{\equ}[1]{(\ref{#1})}
\newcommand{\la}{\lambda}

\newtheorem{theorem}   {Theorem}[section]
\newtheorem{proposition}  {Proposition}[section]

\newtheorem{lemma} {Lemma}[section]

\newtheorem{definition}  {Definition}[section]

\begin{document}

\begin{abstract}
Let $\Omega$ be a bounded  domain in $\RR^N$, $N\ge 2$, with
smooth boundary $\partial\Omega$. We construct positive weak
solutions of the problem $\Delta u + u^p = 0 $ in $\Omega$, which
vanish in suitable trace sense on $\partial\Omega$, but which are
singular at prescribed single points if $p$ is equal or slightly
above $\frac{N+1}{N-1}$.  Similar constructions are carried out
for solutions which are singular on any given embedded submanifold
of $\partial\Omega$ of dimension $0\le k\le N-2$, if $p$ equals or
it is slightly above  $\frac{N-k+1}{N-k-1}$, and even on countable
families of these objects, dense on a given closed set. The role
of this exponent, first discovered by Brezis and Turner
\cite{Bre-Tur} for boundary regularity when $p< \frac{N+1}{N-1}$,
parallels  that of $p=\frac N{N-2}$ for interior singularities.
\end{abstract}

\maketitle

\setcounter{equation}{0}

\section{Introduction and statement of  main results}

Let $\Omega$ be a bounded domain in ${\mathbb R}^N$, with smooth
boundary $\partial \Omega $.
A model of nonlinear elliptic boundary value problem  is the
classical Lane-Emden-Fowler equation,
\begin{eqnarray}
\left\{
\begin{array}{rlllllll} \Delta u + u^{p} & = & 0 \qquad &  \mbox{in} \quad \Omega  \\[2mm]
                         u & > & 0 \qquad &\mbox{in} \quad \Omega \\[2mm]
                         u &  = & 0 \qquad &\mbox{on} \quad \partial\Omega
\end{array}\right.
\label{eq:1}
\end{eqnarray}
where $p>1$. We are interested in finding solutions to this problem
which are smooth in $\Omega$ and equal to $0$ almost everywhere  on
$\partial \Omega$ with respect to surface measure. More precisely,
we want to study solutions to problem \equ{eq:1} that satisfy the
boundary condition in a suitable trace sense, while not necessarily
in a continuous fashion.

\medskip
Following Brezis \& Turner \cite{Bre-Tur} and  Quittner \& Souplet
\cite{Qui-Sou}, we  say that a positive function $u \in {\mathcal
C}^\infty (\Omega)$ is a {\em very weak solution} of problem
$(\ref{eq:1})$ if
$$u,\, u^p\mbox{dist}\, (x,\partial\Omega) \in L^1(\Omega)$$
 and
\[
\int_\Omega  (u \, \Delta v + u^p \, v ) \,  dx =0 \quad\hbox{for
all $v\in {\mathcal C}^{2}(\bar \Omega)$ with $v=0$ on $\partial
\Omega$.}\] From the results in \cite{Bre-Tur,Qui-Sou}, it follows
that if $p$ satisfies the constraint
\begin{equation}
1< p< \frac{N+1}{N-1} \label{expo}\end{equation} then  a very weak
solution $u$ is  actually in $H_0^1(\Omega)$, and it is a weak
solution in the usual variational sense:
\[
u\in H_0^1(\Omega) ,\quad \int_\Omega  (\nabla u \, \nabla  v  - u^p
\, v ) \, dx =0 \quad\hbox{for all $v\in H_0^1(\Omega)$.}\] Elliptic
regularity then yields $u\in {\mathcal C}^2(\bar\Omega)$, so that
$u$ solves \equ{eq:1} classically. As it is well-known, a
constrained minimization procedure involving Sobolev's embedding
implies existence of a weak-variational solution to \equ{eq:1} for
$1<p< \frac{N+2}{N-2}$.  A natural question  is then whether very
weak solutions of \equ{eq:1} are classical within a broader range of
exponents than \equ{expo}. Partially answering this question
negatively, Souplet \cite{Sou} constructed an example of a positive
function $a\in L^\infty(\Omega)$ such that Problem \equ{eq:1}, with
$u^p$ replaced by $a(x)u^p$ for $p> \frac{N+1}{N-1}$, has a very
weak solution which is {\em unbounded}, developing a point
singularity on the boundary.

\medskip
The exponent $p={N+1\over N-1}$ is thus critical in what concerns to
boundary regularity for very weak solutions. The aim of this paper
is to construct solutions to Problem \equ{eq:1} with prescribed
singularities on the boundary. To state an important special case of
our main results we need a definition:

\begin{definition} Let  $u(x)$ be a function defined in $\Omega$ and $x_0\in
\partial\Omega$. We say that
$$
u(x) \to \ell \quad \hbox{as } x\to x_0\ \hbox{ non-tangentially }
$$
if
$$
\lim_{\Gamma_{\alpha} (x_0) \ni x \to x_0 } u(x) =  \ell  \quad
\foral\, \alpha\in [0, \frac \pi 2 ),
$$
where $\Gamma_\alpha (x_0)$ denotes the cone with vertex $\xi_i$,
and angle $\alpha$ with respect to its axis, the inner normal to
$\partial \Omega$ at $x_0$.
\end{definition}
We have the validity of the following result.

\begin{theorem}\label{teo1}
There exists a number $p_N > {N+1\over N-1}$ such that if $p$
satisfies
$$ {N+1\over N-1} \le p < p_N, $$
then the following holds: given points $\xi_1,\xi_2,\ldots, \xi_k\in
\partial\Omega$, there exists a very weak solution $u$ to problem
$\equ{eq:1}$ such that $u\in C^2(\bar \Omega \setminus
\{\xi_1,\ldots,\xi_k\})$ and
$$
u(x) \to +\infty  \quad \hbox{as } x\to \xi_i \hbox{
non-tangentially, } \foral \ i=1,\ldots, k.
$$
\end{theorem}
The study of the behavior near an isolated boundary singularity of any positive solution of (\ref{eq:1}) when when the exponent $p \geq {N + 1 \over N - 1}$ was recently achieved by  Bidaut-V\'eron-Ponce-V\'eron in \cite{BVPV}.

\medskip

\subsection{The parallel with $p=\frac N{N-2}$ and interior
singularities} The role of the exponent $p={N+1\over N-1}$ parallels
that of $p=\frac N{N-2}$ for solutions to problem \equ{eq:1} with
{\em interior singularities}.
 Let us recall that if $u \in L^p(\Omega)$ is a positive
distributional solution of (\ref{eq:1}) and  $1<p <
\frac{N}{{N-2}}$, then $u$ is smooth in $\Omega$. On the other hand,
for $p \geq \frac{N}{{N-2}} $, distributional solutions of
(\ref{eq:1}) with prescribed interior singularities are built in
\cite{Gid-Spr,Pac1,Pac,Chen-Lin,Maz-Pac,Facki1,Facki2}.
Basic cells in those  constructions   are radially symmetric
singular solutions $u=u(|x|)$ for the equation
\begin{equation} \Delta u + u^p = 0 .\label{eq0}\end{equation}
Whenever $p>\frac N{N-2}$, the function
\begin{equation}
u_0(|x| ) =  c_{p,N}\,  |x|^{-\frac 2 {p-1}}\, ,  \qquad c_{p,N} = {
\left [ \frac 2{p-1}( N-2 - \frac 2{p-1})\right ]^{\frac 1{p-1}}},
\label{u0}\end{equation}
 is a explicit singular solution of
\equ{eq0} in $\RR^N\setminus\{0\}$. If, in addition, $\frac N{N-2} <
p < \frac {N+2}{N-2},$ phase plane analysis for the ODE
corresponding to radial solutions of \equ{eq0}, yields existence of
a singular positive solution $u_1$ which connects the behavior of
$u_0$ near the origin with {\em fast decay at infinity},
\begin{equation}
u_1(|x|) = c_{p,N}\,  |x|^{-\frac 2 {p-1}}(1+o(1)) \quad \hbox{as }
x\to 0, \label{u11}\end{equation}
\begin{equation}u_1(|x|) =  |x|^{ -(N-2) }(1+o(1)) \quad \hbox{as } |x|\to +\infty ,
\label{u12}\end{equation}(note that $N-2> \frac 2{p-1}$). The
scalings $u_\la(r) = \la^{2\over p-1} u_1(\la r)$ with $\la
>0$ are then solutions of \equ{eq0} that have the same behavior near the origin
but which become very small as $\la \to 0^+$ on any compact subset
of $\RR^N\setminus \{0\}$. Thus, given  points
$$ \xi_1,\xi_2,\ldots , \xi_k \in \Omega, $$
the function
$$
u_*(x) = \sum_{i=1}^k u_\la (|x-\xi_i|)
$$
constitutes  a ``good approximation'' for small $\la>0$ to a
singular solution of Problem \equ{eq:1}. Linear theory and
perturbation arguments lead to establish
 the presence of an actual solution to \equ{eq:1}   near $u_*$, see \cite{Maz-Pac}. When $p=\frac
N{N-2}$ a similar construction can be carried out, see \cite{Pac}.
Basic cell  $u_1$ corresponds in this case to a  positive radial
solution $u_1$ of equation \equ{eq0} in $B(0,1)$  with
\begin{equation}
u_1(|x|) = c_N |x|^{-(N-2)} \, \log (1/|x| )^{-\frac {N-2} 2}
(1+o(1)) \quad \hbox{as } x\to 0 . \label{u2}\end{equation} In this
case the scalings $u_\la(x) = \la^{N-2 \over 2}u_1(\la x)$  have the
same behavior as $u_1$ at the origin, and they approach zero as $\la
\to 0^+$, uniformly on compact subsets of $\RR^N\setminus\{0\}$.

\subsection{The basic cells: singular solutions on a half-space}

In the construction of the solutions predicted by Theorem \ref{teo1}
we will follow a scheme similar to that described above for interior
singularities. Basic cells will now be positive solutions of
equation \equ{eq0} defined on the half-space,
\[
{\mathbb R}^N_+ : = \{ x = (x_1, \ldots, x_N) \  / \  x_N >0 \}
\]
which vanish on its boundary, with a singularity at the origin. Such
solutions are of course not radial, and ODE analysis does not apply.
Thus, we consider the following two problems:
\begin{eqnarray}
\left\{
\begin{array}{rlllllll} \Delta u + u^{p} & = & 0 \qquad &  \mbox{in} \quad \RR^N_+\setminus\{0\} \\[2mm]
                         u & > & 0 \qquad &\mbox{in} \quad \RR^N_+ \\[2mm]
                         u &  = & 0 \qquad &\mbox{on} \quad \partial
\RR^N_+\setminus\{0\},
\end{array}\right.
\label{eq2}
\end{eqnarray}
for $p> \frac {N+1}{N-1}$, and
\begin{eqnarray}
\left\{
\begin{array}{rlllllll} \Delta u + u^{\frac{N+1}{N-1}} & = & 0 \qquad &  \mbox{in} \quad B_+ \\[2mm]
                         u & > & 0 \qquad &\mbox{in} \quad B_+ \\[2mm]
                         u &  = & 0 \qquad &\mbox{on} \quad
\partial \RR^N_+ \cap \bar B_+\setminus \{0\},
\end{array}\right.
\label{eq3}
\end{eqnarray}
where $B_+ = \RR^N_+ \cap B(0,1) .$

\medskip
Our purpose is to find families of solutions $u_\la$ with analogous
behavior to the radial singular ones  previously described.
 Let us consider first the case $p> \frac {N+1}{N-1}$. The role
of the explicit radial solution $u_0$ in \equ{u0} is now played by
one found by separation of variables: Let us denote by $S^{N-1}_+$
the half sphere
\[
S^{N-1}_+ : = \{ \theta = (\theta_1, \ldots, \theta_N)\in S^{N-1} \
/ \ \theta_N
>0 \}.
\]
Looking for a solution of problem \equ{eq2} of the form
\begin{equation}
u_0(x) =  r^{-\frac 2{p-1}} \phi_p( \theta_N ),\quad r= |x|,\quad
\theta= \frac {x}{|x|}, \label{u00}\end{equation} we arrive at the
problem on the half sphere,
\begin{eqnarray}
\left\{
\begin{array}{rlllllll} \left ( \Delta_{S^{N-1}} + N-1 \right) \, \phi_p  -
\frac{_{p+1}}{^{p-1}} \, \left( N - \frac{_{p+1}}{^{p-1}} \right)
\, \phi_p + \phi_p^{p}  & = & 0 \qquad &  \mbox{in} \quad S^{N-1}_+   \\[2mm]
                         \phi_p & > & 0 \qquad &\mbox{in} \quad S^{N-1}_+ \\[2mm]
                         \phi_p &  = & 0 \qquad &\mbox{on} \quad \partial S^{N-1}_+
.
\end{array}\right.
\label{phip}
\end{eqnarray}
Here $\Delta_{S^{N-1}}$ designates the Laplace-Beltrami operator in
$S^{N-1}$. Since $N-1$ is the first eigenvalue of
$-\Delta_{S^{N-1}}$ under Dirichlet boundary conditions, with
eigenfunction $\theta_N$, in the considered range $N -
\frac{_{p+1}}{^{p-1}}>0$, an application of the mountain pass lemma
yields existence of a solution to this problem, provided that,
additionally, $p$ is subcritical in dimension $N-1$, namely $p<
\frac{N+1}{N-3}$. When $p$ tends from above to $\frac{N+1}{N-1}$,
this solution ceases to exist by uniform vanishing.  Alternatively,
in this regime,  a standard application of Crandall-Rabinowitz local
bifurcation theorem yields that this solution defines a continuous
branch in $p$ with asymptotic behavior \begin{equation}
\phi_p(\theta_N )\, =\, c_N\, (N - \frac{_{p+1}}{^{p-1}})^{\frac
1{p-1}}\, \theta_N\, (1+ o(1)),\quad \hbox{as } p\downarrow
{N+1\over N-1} \, . \label{expansion}\end{equation} Nevertheless,
the function $u_0$ does not suffice  for the construction of
approximate profiles for those of Theorem \ref{teo1} since it is
``too large'' at infinity. We need an analogue of the radial
function $u_1$ in \equ{u11}-\equ{u12}, namely one that  behaves like
$u_0$ near the origin but having {\em fast decay}. A ``connection''
between $u_0$ with  Poisson's kernel $x_N/|x|^N$ does indeed exist
provided that $p$ is sufficiently close to $N+1\over N-1$, as the
following result states.
\begin{proposition} There exists a number $p_{N} > \frac{N+1}{N-1}$, such that for all $\frac{N+1}{N-1} <p<
p_N$, there exists a solution $u_1(x)$
to problem $\equ{eq2}$ such that
\[
u_1(x)  \,= \, |x|^{-\frac 2{p-1}}\, \phi_p ( x_N/|x|)\,(1+ o(1))
\quad \mbox{\em as $x\to 0$},
\]
where  $\phi_p$ solves $\equ{phip}$, and
\[
 u_1(x) \, =\,   |x|^{-N} x_N  \, (1+o(1))\quad \mbox{\em as $|x|\to +\infty$
}.
\]
\label{prop1}
\end{proposition}
This solution  has indeed ``fast decay'' since $N-1 > \frac 2{p-1}$.
Observe then that the scalings $u_\la (x) = \la^{\frac 2{p-1}}
u_1(\la x)$ define a family of solutions to Problem \equ{eq2} which
have a common, $\la$-independent behavior at the origin, but which
vanishes uniformly as $\la\to 0$, on compact subsets of
$\RR^N_+\setminus \{0\}$.

\medskip
When $p=\frac{N+1}{N-1}$ there is no solution to problem \equ{phip}
and thus separation of variables fails. On the other hand, we have
an exact analogue of the radial solutions $u_1$ in \equ{u2}, as
described by the following result.

\begin{proposition}
There exists a  solution $u_1$  of Problem $\equ{eq3}$ such that
\[
u_{1}  (x) =  c_N \, |x|^{-N} \, \log (1/|x|)^{\frac{1-N}{2}} \, x_N
\,(1+ o(1)) \quad\hbox{\em as } x\to 0.
\]
\label{prop2}
\end{proposition}
We observe that in this case the functions $u_\la (x) =
\la^{N-1}u_1(\la x)$ satisfy that $u_\la (x)\to 0$ uniformly as $\la
\to 0^+$ on compact subsets of $\RR^N_+ \setminus \{0\}$.

\bigskip
\subsection{Solutions with prescribed singular set: general statements}
 In reality, the profiles given by the above results can also be
used to approximate solutions to Problem \equ{eq:1} whose singular
set is a $k$ dimensional submanifold of $\partial \Omega$ with $1\le
k\le N-2$. For instance, if $u_1(x')$, $x'\in \RR^{N-k}_+$ is the
solution of \equ{eq2} given by Proposition \ref{prop1} for $p$ close
from above to $\frac{N-k+1}{N-k-1}$, then $\tilde u (x) = u_1(x')$
solves the same problem in $\RR^N_+$, now with singular set given by
a $k$-dimensional subspace. This is the content of the following
result,  more general than Theorem \ref{teo1}, whose analogue for
interior singularities was found in \cite{Pac,Maz-Pac}.

\begin{theorem}
Let $0\le k\le  N-2$ and let $p_{N-k}$ be the number given by
Proposition \ref{prop1} with $N$ replaced by $N-k$. Given $p$ with
$$
\frac{N-k+1}{N-k -1} \le p < p_{N-k}
$$
and a $k$-dimensional submanifold $S$ embedded in $\partial \Omega$,
there exist infinitely many (very) weak solutions to problem
$\equ{eq:1}$ such that
 $u \in {\mathcal C}^{2} (\bar \Omega \setminus S)$, and
$$
u(x) \to +\infty  \quad \hbox{as } x\to x_0 \hbox{ non-tangentially,
} \foral x_0 \in S.
$$
 \label{th:1}
\end{theorem}

When $k=0$, we agree that $S$ is  a finite set of isolated points,
so that Theorem \ref{teo1} is recovered.  In reality, the solutions
found arise as continua, depending on as many real parameters as
number of points lie in $S$.   When $k \geq 1$,  the solutions we
construct are infinite dimensional families. The construction
actually allows much more: For instance, when $p={N+1\over N-1}$,
the number of points of the singular set can be taken to infinity,
to total a dense subset of any given closed set $\mathcal A$ of
$\partial \Omega$, which can be properly called its {\em singular
set}. In fact, since the solutions we are interested in are smooth
in $\Omega$, it is natural to define the singular set of a very weak
solution $u$ of (\ref{eq:1}) as the complement in $\partial \Omega$
of the set of points $x \in
\partial \Omega$ in a neighborhood of which $u$ is smooth. Observe
that, by definition, the singular set of $u$ is a closed subset of
$\partial \Omega$.  We have the validity of the following general
result.

 \begin{theorem}
Let $0\le k\le  N-2$, and $ \frac{N-k+1}{N-k -1} \le p < p_{N-k}. $
Let us consider  a nonempty closed subset $\mathcal A$ of $\partial
\Omega$, which contains a sequence of $k$-dimensional embedded
submanifolds $S_i$, $i\in {\mathbb N}$, which are also disjoint and
satisfy that $S : = \cup_i S_i$ is dense in $\mathcal A$. Then,
there exists a positive very weak solution of Problem $\equ{eq:1}$
whose singular set is exactly $\mathcal A$, and such that
$$
u(x) \to +\infty  \quad \hbox{as } x\to x_0 \, \hbox{
non-tangentially, } \foral x_0 \in S,
$$
and
$$
u(x) \to  0  \quad \hbox{as } x\to x_0 \, \hbox{ non-tangentially, }
\foral x_0 \in \partial\Omega \setminus S.
$$
 \label{th:2}
\end{theorem}

\medskip
This last result and the underlying construction have interesting
consequences: for instance, for $p$ larger than but close enough to
$ \frac{{N+1}}{{N-1}}$, there are infinitely many very weak
solutions of \equ{eq:1} whose singular set is any prescribed closed
subset of $\partial\Omega$, but such that $ u \in W^{1,q}_0(\Omega)$
for any $ 1< q  < N \, \frac{_{p-1}}{^{p+1}} . $ Therefore, even
though $u$ is not identically equal to $0$ at each point of
$\partial \Omega$, we can say that $u=0$ on $\partial \Omega$ in
appropriate sense of traces.

\medskip
The proof of these results relies on two basic ingredients: one is
the construction  of the basic cells of Propositions \ref{prop1} and
\ref{prop2}, which we carry our in \S 2. The other ingredient is the
analysis of invertibility of Laplace's operator for right hand sides
that involve singular behavior near a point or an embedded manifold
of the boundary. After this analysis, which is carried out in \S 3,
the proof of Theorem \ref{th:1} then follows from a fixed point
argument. The result of Theorem \ref{th:2} is a consequence of an
inductive construction taken to the limit under suitable control.

\bigskip
\setcounter{equation}{0}
\section{The half-space case: proofs of Propositions \ref{prop1} and \ref{prop2}}

It is natural and convenient to look for solutions of (\ref{eq2}) or
\equ{eq3} of the form
\[
u(x) = |x|^{-\frac{2}{p-1}}  \, \phi ( -\log |x|, x/|x| ) ,
\]
so that the equation $\Delta u + u^p=0$ reads in terms of
$\phi(t,\theta)$, $t \in {\mathbb R}$, $\theta \in S_+^{N-1}$, as
\begin{equation}
\partial_t^2 \phi  - \left( N - 2 \frac{_{p+1}}{^{p-1}}\right) \, \partial_t \phi -
\frac{_{p+1}}{^{p-1}} \, \left( N  - \frac{_{p+1}}{^{p-1}} \right)
\, \phi + \left( \Delta_{S^{N-1}} + N-1 \right) \, \phi + \phi^{p} =
0 .\label{eq5}
\end{equation}
\medskip

\subsection{Proof of Proposition~\ref{prop2}}

When $p = \frac{{N+1}}{{N-1}}$, in the language of equation
(\ref{eq5}), problem \equ{eq3} becomes

\begin{eqnarray}
\left\{
\begin{array}{rlllllll} \partial_t^2 \phi +  N  \, \partial_t \phi + \left(
\Delta_{S^{N-1}} + N-1 \right) \, \phi + \phi^{\frac{N+1}{N-1}} & = & 0 \quad &  \mbox{in} \quad (t_*,\infty)\times S^{N-1}_+   \\[2mm]
                         \phi & > & 0 \quad &\mbox{in} \quad (t_*,\infty)\times S^{N-1}_+ \\[2mm]
                         \phi &  = & 0 \quad &\mbox{on} \quad  (t_*,\infty)\times
                         \partial S^{N-1}_+\,
                         .
\end{array}\right.
\label{eq33}
\end{eqnarray}
We allow here $t_*> 0$ to be a parameter, which we will choose later
to be large. To get a solution of problem \equ{eq3} we actually need
$t_* =0$, but this is simply achieved by a translation of $\phi$ in
the $t$-variable.

\medskip
 The idea is now to look for a solution of this
equation of the form
\begin{equation}
\phi (t,\theta) = a_N \, t^{-b_N} \, \ee(\theta) +  \psi (t,\theta)
, \label{psi}\end{equation}  where $a_N, b_N$ are positive constants
to be fixed below, and $\ee$ denotes the eigenfunction of $-
\Delta_{S^{N-1}_+}$ associated to the eigenvalue $N-1$ and
normalized so that its $L^2$-norm is equal to $1$. Explicitly,
$$
\ee(\theta) = \frac {\theta_N}  {\left (\int_{S_+^{N-1}} \theta_N^2
\, d\sigma \right )^{1/2}}.$$  When substituting the function $\phi
=  a_N t^{-b_N} \, \ee(\theta)$ as an approximation for a solution
of  equation \equ{eq33}, we see that for large $t$, the main order
term in the error created is the function
\[
E(t,\theta) : = \,  N  a_N  b_N \, t^{-b-1} \, \ee (\theta) - |a_N\,
t^{-b_N} \, \ee (\theta) |^{\frac{N+1}{N-1}} .
\]
We make the following choice for the numbers $a_N$ and $b_N$:
\[
b_N = \frac{_{N-1}}{^2}, \qquad a_N =  \left [\, \frac 2{{N(N-1)}}
\int_{S^{N-1}_+} \ee^{\frac{2N}{N-1}} \, d\sigma \,\right ]^{-\frac
{N-1}2}
 \,        .
\]
This election achieves the $L^2$-orthogonality of $E$ to $\ee$ for
all $t$, namely
$$
\int_{S^{N-1}_+} E(t,\cdot)\, \ee \, d\sigma = 0\quad\hbox{for all }
t> t_* .
$$
We fix  these  values in what follows. In terms of $\psi$ in
\equ{psi}, equation \equ{eq33} now reads
\[
 \partial_t^2 \psi +  N  \, \partial_t \psi + \left(
\Delta_{S^{N-1}} + N-1 \right) \, \psi +  b_N \, (b_N+1) \, a_N \,
t^{-b_N-2} \, \ee  -\qquad\qquad
\]
\[
\qquad\quad \qquad N  b_N a_N \, t^{-b_N-1} \, \ee  + |a_N\,
t^{-b_N} \, \ee + \psi|^{\frac{N+1}{N-1}} =0,
\]
$$
                         \psi   = 0 \quad \mbox{on} \quad  (t_*,\infty)\times
                         \partial S^{N-1}_+ . $$
We further decompose \begin{equation} \psi (t,\theta) =  f_2(t) \,
\ee(\theta) + \psi_1 (t, \theta )\label{psi1}\end{equation} where
$\psi_1(t, \theta)$ satisfies
$$
\int_{S^{N-1}_+} \psi_1(t,\cdot)\, \ee \, d\sigma = 0\quad\hbox{for
all } t> t_* .
$$
The equation we have to solve then reduces to the coupled system in
$(\psi_1,f_2)$ given by
\begin{equation}
\left\{ \begin{array}{rlllllll} \left( \partial_t^2 +  N  \,
\partial_t  + \left( \Delta_{S^{N-1}} + N-1 \right) \right) \,
\psi_1 & = &  N_1(\psi_1,f_2) \\[3mm]
( \partial_t^2 +  N  \, \partial_t  + \frac{N(N-1)}{2} \,
\frac{1}{t} ) \, f_2 & = & N_2(\psi_1,f_2) ,\\[3mm]
 \psi_1  & = &0 \quad \mbox{on} \quad  (t_*,\infty)\times
                         \partial S^{N-1}_+ ,
\end{array}
\right. \label{eq:2.9}
\end{equation}
where
\begin{equation}
\begin{array}{rlllllll} N_1(\psi_1,f_2) & = &  ( \frac{N(N-1)}{2} \, a_N \, \ee -
a_N^{\frac{N+1}{N-1}}\, \ee^{\frac{N+1}{N-1}} ) \,t^{-
\frac{N+1}{2}} \\[3mm]
& - & \Pi^{\perp} \left( |a_N\, t^{\frac{1-N}{2}} \, \ee + \psi_1
+f_2\ee  |^{\frac{N+1}{N-1}} -  |a_N\, t^{\frac{1-N}{2}} \, \ee
|^{\frac{N+1}{N-1}} \right)\, ,\\[3mm]
N_2(\psi_1,f_2) & = & \frac{N^2-1}{4} \, a_N \,
t^{-\frac{N+3}{2}} - \\[3mm]&&
\displaystyle \int_{S^{N-1}_+} (\,  |a_N\, t^{\frac{1-N}{2}} \, \ee
+ \psi_1 +f_2\ee |^{\frac{N+1}{N-1}}  - |a_N\, t^{\frac{1-N}{2}} \,
\ee
|^{\frac{N+1}{N-1}} \\[3mm]&&
\qquad\quad - \frac{{N+1}}{{N-1}} \, a_N^{\frac{2}{N-1}} \,
\ee^{\frac{N+1}{N-1}} \, \frac{1}{t}  f_2 \, ) \, \ee \, d\sigma \,
.
\end{array} \label{N}
\end{equation}
Here $\Pi^\perp$ denotes the $L^2$-orthogonal projection over the
orthogonal complement to $\ee$, namely
$$
\Pi^\perp (h)  = h(t,\theta ) -  \ee(\theta) \int_{S^{N-1}_+}
h(t,\cdot)\, \ee \, d\sigma ,
$$
and $\psi$ is given by \equ{psi1}. The logic in the resolution of
problem \equ{eq:2.9} is simple: we look for a solution
$(\psi_1,f_2)$ which is  small compared with $ t^{\frac{1-N}{2}} \,
\ee$. We will construct inverses to the linear operators defined by
the left hand sides of the equations in \equ{eq:2.9} with suitable
bounds that allow, for sufficiently large $t_*$, the resolution of
the system via contraction mapping principle. Observe that so far we
have not imposed boundary conditions at $t=t_*$. We will invert the
linear operator in $\psi_1$, for right hand sides $L^2$-orthogonal
to $\ee$ for all $t$, imposing Dirichlet boundary condition at
$t=t_*$. The choice of inverse for the ODE operator in $f_2$ will be
basically explicit, and will not require imposing boundary
conditions. The natural environment to carry out these inversions is
$L^\infty$-weighted spaces. In the next two lemmas we construct
these inverses. Thus we consider the linear problems
\begin{equation}
\left\{ \begin{array}{rlllllll} \left( \partial_t^2 +  N  \,
\partial_t  + \left( \Delta_{S^{N-1}} + N-1 \right) \right) \,
\psi & = &   h  \quad \mbox{in} \quad  (t_*,\infty)\times
                               S^{N-1}_+ ,\\[3mm]
 \psi  & = &0 \quad \mbox{on} \quad  \partial\, (\, (t_*,\infty)\times
                          S^{N-1} ) \\ [3mm]
 \int_{S^{N-1}_+} \psi(t,\cdot)\, \ee \, d\sigma & =&0\quad\hbox{for all } t> t_* ,
\end{array}
\right. \label{eq:2.5}
\end{equation}
for $h$ such that \begin{equation}\int_{S^{N-1}_+} h(t,\cdot)\, \ee
\, d\sigma =0\quad\hbox{for all } t> t_* , \label{orto}
\end{equation} and
\begin{equation}
\left( \partial_t^2 +  N \, \partial_t + \frac{_{N(N-1)}}{^2} \,
\frac{_1}{^t} \right) \, f = g \quad \mbox{in} \quad
(t_*,\infty).\label{eq:2.11}
\end{equation}
We have the validity of the following results.

\medskip
\begin{lemma}
\label{lemaT} There exists a  constant $c>0$ such that the following
holds:  Given $\sigma \geq0$, there is a  $t_\sigma
>0$, with $t_\sigma >0$ if $\sigma >0$ and $t_\sigma =0$ if $\sigma =0$,
 such that, for all $t_* \geq t_\sigma$, and all $h \in
{\mathcal C}^0 ((t_*, \infty)\times S^{N-1}_+)$ that satisfies
$\equ{orto}$ and $ t^{\sigma}h \in L^\infty ((t_*, \infty)\times
S^{N-1}_+) $, there exists a solution $\psi= T_1(h)$ of problem
$\equ{eq:2.5}$, which defines a linear operator of $h$ and satisfies
the estimate
\[
\| t^{\sigma}  \, \psi \|_{L^\infty} + \| t^{ \sigma}  \,
\nabla_\theta \psi \|_{L^\infty} \, \leq \, c \, \| t^{\sigma} \, h
\|_{L^\infty } .
\]
\end{lemma}

\begin{lemma}
\label{lemaTT} Given $\sigma  > \frac{_{N-1}}{^2}$, there exist
numbers $t_\sigma ,\, c_\sigma >0$ such that for all $t_*
>t_\sigma$ and all $g \in {\mathcal C}^0 ((t_* , \infty))$
satisfying $t^{\sigma} g  \in L^\infty ((t_* , \infty))$, there
exists a solution $f=T_2(g)$ of equation $\equ{eq:2.11}$, which
defines a linear operator of $g$ and satisfies the estimate

\[
\| t^{\sigma}  \, f \|_{L^\infty} \leq \, c_\sigma \, \|t^{1+\sigma}
\, g \|_{L^\infty } .
\]
\end{lemma}

Before proceeding into the proofs of these lemmas, let us conclude
the result.

\bigskip
\noindent{\bf Conclusion of the proof of Proposition \ref{prop2}.}
Let us fix in the above lemmas any number $\sigma$ such that
$$  \frac{_{N-1}}{^2} <\sigma <  \frac{_{N+1}}{^2}$$ and $t_*>
t_\sigma$. We obtain a solution  of  problem (\ref{eq:2.9}) if
$(\psi_1,f_2)$ solves the  fixed point problem
\begin{equation}
(\psi_1, f_2) \,= \, {\mathcal M}(\psi_1, f_2) := (T_1(\,
N_1(\psi_1,f_2)\,)\,,\, T_2(\, N_2(\psi_1,f_2)\,)) , \quad
\label{fp}\end{equation}
 in the space
of functions
\[
(\psi, f) \in  {\mathcal C}^0 ([t_* , \infty) \times S^{N-1}_+)
\times {\mathcal C}^0 ([t_*, \infty))
\]
for which the norm
$$
\|(\psi,f)\|_\mu =
 \| t^{\sigma} \,
\psi\|_{\infty} +  \mu \| t^{\sigma} f\|_{\infty}$$ is finite. Here
$\mu<1$ is a positive number which we will fix later. $T_1$, $T_2$
are the operators predicted by Lemmas \ref{lemaT} and \ref{lemaTT}.
It is directly checked that we have the pointwise estimates
\begin{equation}
\begin{array}{rlllllll}
|N_1(\psi_1,f_2)| &\le & A[\, t^{ -\frac {N+1}2 } +  t^{-1}
|\psi_1|+ t^{-1}|f_2| \, ]\, ,\\ [3mm]
 |N_2(\psi_1,f_2)| &\le &  A[\,
t^{ -\frac {N+3}2 } +  t^{-1} |\psi_1|+ t^{N-3 \over 2} |f_2|^2 +
|f_2|^{N+1\over N-1}  \, ]\, ,
\end{array}
\end{equation}
where $A$ depends only on $N$, whenever  $ \|(\psi,f)\|_\mu \le
\mu$. It follows that,
\begin{equation}\begin{array}{rllllll}\| t^{\sigma} N_1(\psi_1,f_2)\|_\infty  &\le & A
\,[\, t_*^{\sigma -\frac {N+1}2 } +  t^{-1}_* \|t^\sigma
\psi_1\|_{\infty} +  t^{-1}_* \|t^\sigma f_2\|_{\infty}\, ]\, , \\
[3mm] \| t^{1+\sigma} N_2(\psi_1,f_2)\|_\infty  &\le & A \,[\,
t_*^{\sigma -\frac {N+1}2 } +    \|t^\sigma \psi_1\|_\infty  +
(t_*^{{N-1 \over 2} -\sigma } + t_*^{1 -\frac {2\sigma}{N-1}} )
\|t^{\sigma}f_2\|_\infty \, ]\, .
\end{array}\end{equation}
These estimates, together with Lemmas \ref{lemaT} and \ref{lemaTT},
yield that if $\mu$ is chosen sufficiently small, depending only on
$\sigma$ and $N$, and  $t^*$ is  taken sufficiently large, then the
operator $\mathcal M$ applies the ball $\|(\psi,f)\|_\mu \le \mu$
into itself. A similar estimates shows that, also, $\mathcal M$ is a
contraction mapping with this norm inside this region. Hence there
is a fixed point $(\psi_1, f_2)$ in this ball. The solution obtained
this way renders the function
$$
\phi (t,\theta) = a_N \, t^{- \frac{N-1}2} \, \ee(\theta) + f_2(t)
\, \ee(\theta)+ \psi_1 (t, \theta )
$$
positive in $(t_*,+\infty)\times S^{N-1}_+$, and it is then a
solution of problem \equ{eq33}.
 This completes the proof of Proposition~\ref{prop2}. \qed

\bigskip
Next we carry out the proofs of the lemmas.

\bigskip
\noindent\noindent{\bf Proof of Lemma \ref{lemaT}~.}  Let us
consider first the case $\sigma =0$, so that  $h$ is bounded.  With
no loss of generality, we also assume $t_*=0$. We see then that
problem \equ{eq:2.5} has at most one bounded solution. This can be
shown for instance expanding a bounded solution of the equation with
$h=0$  in eigenfunctions of the Laplace-Beltrami operator with zero
boundary conditions on $S_+^{N-1}$. The coefficients in this
expansion will be functions of $t$ which correspond to bounded
solution of certain homogeneous ODE's which only have the zero
solution. Thus, we only have to prove existence. To do so, let us
consider, for any given number $t_2>0 $, the problem
\begin{equation}\begin{array}{rllllll}
(\partial_t^2 +  N \, \partial_t + \left( \Delta_{S^{N-1}} + N-1
\right) ) \, \psi &= & h \quad \hbox{in }(0, t_2) \times S^{N-1}_+
\\ [3mm]
\psi &= & 0 \quad \hbox{on } \partial\,(\; (0, t_2) \times
S^{N-1}_+).
\end{array}\label{t2}\end{equation}
This problem is uniquely solvable since it is just a rephrasing of a
Dirichlet problem for the Laplacian in a half-annular region. Let us
denote by $\psi = \psi_{t_2}$ its unique solution. Since, by
assumption, $h (t, \cdot)$ is $L^2$-orthogonal to $\ee$ for all $t
\in (0 , t_2)$, so is $\psi$.

\medskip
It suffices to check that there exists a constant $c >0$ independent
of $t_2 \geq 1$ such that
\begin{equation}
\| \psi \|_{L^\infty( [0,t_2] \times S^{N-1}_+)} \leq \, c \, \| h
\|_{L^\infty ([0,t_2] \times S^{N-1}_+) } .
\label{equation}\end{equation} Indeed, assuming this estimate is
already proven, we use elliptic estimates together with Ascoli's
theorem to show that, as $t_2$ tends to $\infty$,  the sequence of
functions $\psi_{t_2}$ converges uniformly to a function $\psi$
solution of (\ref{eq:2.5}) which  satisfies
\[
\| \psi \|_{L^\infty( [0,\infty) \times S^{N-1}_+)} \leq \, c \, \|
h \|_{L^\infty ([ 0,\infty ) \times S^{N-1}_+) } .
\]
 Elliptic estimates
then imply that
\begin{equation}
\|\nabla  \psi \|_{L^\infty( [0,\infty) \times S^{N-1}_+)} + \| \psi
\|_{L^\infty( [0,\infty) \times S^{N-1}_+)} \leq \, c_0 \, \| h
\|_{L^\infty ([0,\infty ) \times S^{N-1}_+) } . \label{eq:2.55}
\end{equation}
The orthogonality conditions on $\psi$ pass certainly to the limit,
and existence of a solution with the desired properties thus
follows. It remains to prove the uniform estimate \equ{equation}. We
argue by contradiction. Since the result is certainly true when
$t_2$ remains bounded, we assume that there exists a sequence
$t_2=t_{2,i}$ tending to $\infty$,  functions $h=h_i$ and $\psi_i$
corresponding solutions to problem \equ{t2} for which
\[
\| \psi_i \|_{L^\infty( [0, t_{2,i}] \times S^{N-1}_+)} = 1 \qquad
\mbox{and} \qquad \lim_{i\rightarrow \infty} \| h_i \|_{L^\infty
([0, t_{2,i}] \times S^{N-1}_+) } =0.
\]
We choose $t_i \in (0, t_{2,i})$ where $\| \psi_i \|_{L^\infty(
[t_{1,i}, t_{2,i}] \times S^{N-1}_+)}$ is achieved and define
\[
\tilde \psi_i (t,\theta) = \psi_i (t+t_i,\theta)
\]
Using elliptic estimates together with Ascoli's theorem, we can
extract from  $(\tilde \psi_i)_i$ some subsequence which converges
uniformly on compact sets  to $\tilde \psi$, a bounded solution of
\begin{equation}
(\partial_t^2 +  N \, \partial_t + \left( \Delta_{S^{N-1}} + N-1
\right) ) \, \tilde \psi = 0 \label{eq:2.6}
\end{equation} which is
either defined on $[0, \infty)\times S^{N-1}_+$, on $(-\infty,
0]\times S^{N-1}_+$ or on $(-\infty , \infty)\times S^{N-1}_+$.
Furthermore,
\begin{equation}
\| \tilde \psi \|_{L^\infty} =1
 \label{eq:2.7}\end{equation}
with $\tilde \psi$ having $0$ boundary data. Furthermore $\tilde
\psi (t, \cdot)$ is $L^2$-orthogonal to $\ee$, for all $t$.
Eigenfunction decomposition of $\tilde \psi(t, \cdot)$ for the
Laplace-Beltrami operator yields that there is non nontrivial
bounded solution of (\ref{eq:2.6}) and this contradicts
(\ref{eq:2.7}). This completes the proof of the uniform estimate,
and thus existence of a unique bounded solution of \equ{eq:2.5} with
the desired estimate follows. This solution of course defines a
linear operator on bounded $h$.

\medskip
To establish the result for $\sigma>0$ and $t_*>0$ sufficiently
large, let us write
\[
h = t^{-\sigma}\, \tilde h \qquad \mbox{and} \qquad \psi =
t^{-\sigma}\, \tilde \psi
\]
so that  $\tilde h$ is bounded. (\ref{eq:2.5}) reduces to
\begin{equation}
(\partial_t^2 +  N \, \partial_t + \left( \Delta_{S^{N-1}} + N-1
\right) ) \, \tilde \psi  +\left( \frac{_{\sigma (\sigma
+1)}}{^{t^2}} - \frac{_{N\sigma}}{^t} \right) \, \tilde \psi -
\frac{_{2\sigma}}{^t} \, \partial_t \tilde \psi = \tilde h
\label{eq:2.66} \end{equation} We can estimate
\[
\begin{array}{rlllllll}
\displaystyle \| \left( \frac{_{\sigma(\sigma +1)}}{^{t^2}} -
\frac{_{N \sigma}}{^t} \right) \, \tilde \psi -
\frac{_{2\sigma}}{^t} \,
\partial_t \tilde \psi \|_{L^\infty ([ t_*,\infty ) \times
S^{N-1}_+) } \leq  \qquad \qquad  \qquad \qquad \qquad \qquad \\[3mm]
\qquad \qquad  \displaystyle  \mu \,  \left( \|\nabla \tilde \psi
\|_{L^\infty( [t_*,\infty) \times S^{N-1}_+)} + \| \tilde \psi
\|_{L^\infty( [t_*,\infty) \times S^{N-1}_+)} \right)
\end{array}
\]
where $\mu$ can be taken as small as we wish, after choosing $t_*
\geq t_\sigma$ with $t_\sigma$  large enough.  The resolution of
(\ref{eq:2.66}) with the desired bound then follows from that of
(\ref{eq:2.5}) with $\sigma =0$ together with a direct linear
perturbation argument. This finishes the proof. \qed

\bigskip
\noindent\noindent{\bf Proof of Lemma \ref{lemaTT}~.}  Observe that
a right inverse for the operator $
\partial_t^2 + N \, \partial_t
$ on $[t_*, \infty)$ is given by
\[
G(g) (t) = - \int_t^\infty e^{-N \zeta} \int_{t_*}^\zeta \, e^{Ns}
\, g (s)\, ds \, dt
\]
One checks that
\[
\| t^{\sigma} \, G(g ) \|_{L^\infty ((t_*, +\infty))} \leq
\frac{_1}{^{N \, |\sigma| }} \, \left(  1- \frac{_{(\sigma+1)}}{^{N
t_*}} \right)^{-1} \, \| t^{1+\sigma} g \|_{L^\infty ((t_*,
+\infty))}
\]
provided $ Nt_*-1-\sigma >0$. This follows at once from the
computation
\[
\begin{array}{rlllll}
\displaystyle \int_{t_*}^t e^{Ns} \, s^{-\sigma-1}\, ds  & = &
\displaystyle \left[ \frac{_1}{^N} \, e^{Ns} \, s^{-\sigma-1}
\right]_{t_*}^t + \frac{_{\sigma+1}}{^N} \, \int_{t_*}^t e^{Ns} \,
s^{-\sigma -2}\, ds
\\[3mm]
 & \leq  &  \displaystyle \frac{_1}{^N} \, e^{Nt} \, t^{-\sigma-1} +
\frac{_{1 +\sigma}}{^N t_*} \, \int_{t_*}^t e^{Ns} \, s^{-\sigma
-1}\, ds
\end{array}
\]
and hence
\[
 \left( 1 - \frac{_{\sigma+1}}{^N t_*} \right)  \int_{t_*}^t e^{Ns}
 \, s^{-\sigma-1}\, ds \leq  \frac{_1}{^N} \, e^{Nt} \,
 t^{-\sigma-1}.
\]
The result of the lemma then follows from a simple linear
perturbation argument, provided $t_* >0$ is chosen so that
\[
\frac{_{N(N-1)}}{^2} \, \| t^{\sigma} \, \frac{1}{t} \, G( g)
\|_{L^\infty ((t_*, +\infty))} \leq \frac{_1}{^2} \, \| t^{1+\sigma}
g \|_{L^\infty ((t_*, +\infty))},
\]
and the result is concluded. \hfill $\Box$

\medskip

\subsection{Proof of Proposition~\ref{prop1}}

Recall that, when $p \in (\frac{_{N+1}}{^{N-1}},
\frac{_{N+1}}{^{N-3}})$ the Mountain Pass Lemma yields the existence
of $\phi_p$, a nontrivial positive solution of (\ref{phip}). This
solution then induces a solution
\[
 u_{p} ( x ) : = |x|^{-\frac{2}{p-1}} \, \phi_p (x/|x|),
\]
of  Problem \equ{eq2}, for which this time we emphasize its
dependence on $p$. We have to show that there exists a solution of
(\ref{eq2}) which is asymptotic to $ u_{p}$ near $0$ and it is
asymptotic to
\[
u_\infty (x) : = |x|^{-N} \,x_N
\]
at infinity. Let us consider a smooth cut-off function $\chi$ which
is  equal to $0$ in $B_1(0)$ and it is identically equal to $1$ in
${\mathbb R}^N \setminus B_2(0)$. We will consider the function $
u_{p}(1-\chi) $ as a first approximation for the solution we are
looking for. Since, we recall, $\phi_p$ approaches zero uniformly as
$p\downarrow {N+1\over N-1}$, then the same is true for $u_{0,p}$
away from the origin. The result of the proposition  relies on a
perturbation procedure, and this is the reason why we can only show
the  for exponents $p$ close to $\frac{_N}{^{N-1}}$.  To carry out
this scheme, we shall build a right inverse for the Laplacian
relative to the following
 doubly weighted space:

\begin{definition}
Given $\delta, \delta' \in {\mathbb R}$, the space
$L^\infty_{\delta, \delta'} ({\mathbb R}^N_+)$ is defined to be the
space of functions $u \in L^\infty_{loc} ({\mathbb R}^N_+)$ for
which the following norm
\[
\| u \|_{L^\infty_{\delta, \delta'} ({\mathbb R}^N_+)} = \|
|x|^{-\delta} \, u\|_{L^\infty (B_+(1))} + \| |x|^{-{\delta'}} \,
u\|_{L^\infty ({\mathbb R}^N_+ - B_+(1))}
\]
is finite.
\end{definition}
Hence $\delta$ controls the behavior of the function near $0$ and
$\delta'$ the behavior of the function near infinity. Let us
consider the problem

\begin{equation}\begin{array}{rlllll}
 \, \Delta \, u &=& |x|^{-2}f \quad \hbox{in } \RR^N_+ \\ [3mm]
               u &=& 0 \quad \hbox{on } \partial\RR^N_+
               \setminus\{0\}.
\end{array}\label{problem}\end{equation}

\medskip
We have the validity of the following result.

\begin{lemma}
Let $\delta \in (1-N, 1)$ and $\delta' \in (-N, 1-N)$ be given.
There is a constant $c>0$  such that for each  $f \in
L^\infty_{\delta, \delta'} ({\mathbb R}^N_+)$,  there exists a
solution $u=G(f)$ of problem $\equ{problem}$, which defines a linear
operator in $f$ and can be decomposed  as $$u = \tilde u + a\chi
u_\infty$$ with
$$
|a| + \|\tilde u\|_{L^\infty_{\delta, \delta'}} \le
c\|f\|_{L^\infty_{\delta, \delta'}} .
$$
 \label{pr:3}
\end{lemma}
\noindent{\bf Proof.} Let us observe that  $\delta \, ( N - 2 +
\delta ) < N-1$ precisely when $\delta \in (1-N, 1)$. Therefore, we
can define $\varphi_* = \varphi_{N, \delta}$ to be the unique,
positive solution of
\begin{equation}\begin{array}{rlllll}
 \,- \left( \Delta_{S^{N-1}}  + \delta  \, (\delta + N -2)
\right) \, \varphi_*&=& 1 \quad \hbox{in } S^{N-1}_+ \\
[3mm]
               \varphi_* &=& 0 \quad \hbox{on } \partial
               S^{N-1}_+.
\end{array}\nonumber\end{equation}
A direct computation shows that
\begin{equation}
- |x|^{2} \, \Delta_{{\mathbb R}^N} ( |x|^{\delta} \, \varphi
_*(\theta)) = c |x|^{\delta }. \label{eq:2.12}
\end{equation}
Assume that $f\in L^\infty_{\delta, \delta'} ({\mathbb R}^N_+)$.
Given $r_1 < 1 < r_2$, we can solve the equation $|x|^2 \, \Delta u
= f$ in ${\mathbb R}^N_+ \cap (B_{r_2} - B_{r_1})$, with $0$
boundary conditions. We use the function $x \longmapsto \varphi_*
(x) \, |x|^\delta$ as a barrier to prove the pointwise estimate
\[
|u| \leq c \,  \| f\|_{L^\infty_{\delta, \delta'} ({\mathbb R}^N_+)}
\, |x|^{\delta}
\]
in ${\mathbb R}^N_+ \cap (B_{r_2} - B_{r_1})$. Furthermore, given
$\bar \delta \in (1-N, 1)$ we can use the function $x \longmapsto
\varphi_*(x) \, |x|^{\bar \delta}$ as a barrier to prove the
estimate
\[
|u| \leq c \,  \| f\|_{L^\infty_{\delta, \delta'} ({\mathbb R}^N_+)}
\, |x|^{\bar \delta}
\]
in ${\mathbb R}^N_+ \cap (B_{r_2} - B_{r_1})$.

\medskip

We use elliptic regularity theory as well as Ascoli's theorem to
pass to the limit as $r_1$ tends to $0$ and $r_2$ tends to $\infty$.
We obtain a solution $u$ of problem \equ{problem} which satisfies
the pointwise estimates
\[
|u| \leq c \,  \| f\|_{L^\infty_{\delta, \delta'} ({\mathbb R}^N_+)}
\, |x|^{\delta}
\]
in $({\mathbb R}^N_+ \cap B_1 )-\{0\}$, and
\[
|u| \leq c \,  \| f\|_{L^\infty_{\delta, \delta'} ({\mathbb R}^N_+)}
\, |x|^{\bar \delta}
\]
in ${\mathbb R}^N_+ \cap ({\mathbb R}^N - B_1)$. Finally, the
decomposition of the solution $u$ at infinity into
\[
 u =  \tilde u + a \, u_\infty
\]
where the function $\tilde u$ satisfies
\[
|\tilde u| \leq c \,  \| f\|_{L^\infty_{\delta, \delta'} ({\mathbb
R}^N_+)} \, |x|^{\delta'}
\]
follows easily from Green's representation formula. Moreover, we can
directly compute the value of $a$. Indeed, integration of the
equation over $B_r^+ : = {\mathbb R}^N_+ \cap B_r$ yields, for $r$
large enough,
\[
\int_{B_r^+} f \, |x|^{-2} \, dx  = \int_{\partial B_r^+} \,
\partial_r  u \, d\sigma = \int_{S^{N-1}_+} \,
(\partial_r  \tilde u ) (r \,\theta )  \, r^{N-1} \, d\sigma - a \,
(N-1) \, \int_{S^{N-1}_+ } \,\theta_N \, d\sigma
\]
Passing to the limit as $r$ tends to $\infty$, we  obtain the
identity
\begin{equation}
a \, (N-1) \, \int_{S^{N-1}_+}\theta_N \, d\sigma  =  -
\int_{{\mathbb R}^N_+} f \, |x|^{-2} \, dx ,
\label{formulaa}\end{equation} and the proof is concluded. \qed

\bigskip
\noindent{\bf Conclusion of the proof of Proposition \ref{prop1}.}
To find a solution of problem \equ{eq2}, we write
$$
u = (1-\chi) \bar u_p + v ,
$$
where $\chi$ is a smooth cut-off function which is equal to $0$ in
$B_1$ and identically equal to $0$ in ${\mathbb R}^N - B_2$. Let us
fix numbers $\delta \in (1-N, 1)$, $\delta' \in (-N, 1-N)$ and let
$G$ be the operator defined in Lemma \ref{pr:3}. Then, we obtain a
solution with the required properties if $v$ solves the fixed point
problem
\begin{equation}
v  = - G  \, \left( |x|^{2} ( \Delta (1-\chi) \, \bar u_{p} + |
(1-\chi) \, \bar u_{p} + v|^{p} )\right) \label{fp1}\end{equation}
in the space $L^\infty_{\delta, \delta'} (({\mathbb R}^N_+-\{0\})
\oplus \,\mbox{Span}\, \{ \chi \, u_\infty \})$, and $(1-\chi) u_p +
v>0$.

\medskip
Let us observe that  there exists a constant $c_0 = c(N, \delta,
\delta')
>0$ such that
\begin{equation}
\| |x|^{2}  \, ( \Delta  ((1-\chi) \,\bar  u_{p}) + ((1-\chi) \,
\bar u_{p})^{p} ) \|_{L^\infty_{\delta , \delta'} ({\mathbb R}^N_+
-\{0\} )} \leq c_0 \, \| \phi_p \|_{{\mathcal C}^{2} (S^{N-1}_+)} .
\label{es1}\end{equation} Indeed, since $\Delta \bar u_{p} + \bar
u_{p}^{p} =0$ we have
$$
\Delta  ((1-\chi) \, \bar u_{p}) + ((1-\chi) \,\bar  u_{p})^{p}  = -
\Delta \chi \, \bar u_{p} - 2 \nabla \chi\cdot \nabla \bar u_{p} +
(1-\chi - (1-\chi)^p) \, \Delta \bar u_{p}
$$
 and the estimate
follows at once.

\medskip
On the other hand, if we  assume that $p$ is sufficiently close to
$N+1\over N-1$ from above, we have that
\[
\| \phi_p \|_{{\mathcal C}^{2} (S^{N-1}_+)}\leq 1 ,
\]
and also
\[
\delta > - \frac{_2}{^{p-1}} \qquad \mbox{and} \qquad \delta ' > p
\, (1-N)+ 2 .
\]
Under these constraints, it is not hard to check the existence of a
constant $c = c(N, \delta , \delta' )
>0$ such that
\begin{equation}
\begin{array}{lllll}
\||x|^{2} ( |(1-\chi) \, \bar u_{p} + v_2|^{p} - |(1-\chi) \, \bar
u_{p} + v_1|^{p} )\|_{L^\infty_{\delta , \delta'} ({\mathbb R}^N_+
-\{ 0 \})}
\\[3mm]
\qquad \qquad \qquad \leq c \,  \| \phi_p \|_{{\mathcal C}^{2}
(S^{N-1}_+)}  \,  \|v_2 - v_1\|_{L^\infty_{\delta, \delta'}
({\mathbb R}^N_+ - \{0\} ) \oplus \, \mbox{Span} \, \{ \chi \,
u_\infty \}}
\end{array}
\label{es2}\end{equation}  for all $v_2, v_1 \in L^\infty_{\delta ,
\delta'} ({\mathbb R}^N_+ - \{0\}) \oplus \mbox{Span} \{\chi \,
u_\infty \}$ satisfying
\[
\| v_i\|_{L^\infty_{\delta, \delta'} (({\mathbb R}^N_+-\{0\})
\oplus\, \mbox{Span}\, \{ \chi \, u_\infty \})} \leq \, 2 \, c_0 \,
\| \phi_p \|_{{\mathcal C}^{2} (S^{N-1}_+)} .
\]

Using  estimates \equ{es1}, \equ{es2} and Lemma \ref{pr:3}, the
existence of a solution to the fixed point  problem  \equ{fp1} can
then be obtained by contraction mapping principle in the ball of
radius $2 \, c_0 \, \| \phi_p \|_{{\mathcal C}^{2} (S^{N-1}_+)}$ in
the space $L^\infty_{\delta, \delta'} (({\mathbb R}^N_+-\{0\})
\oplus \,\mbox{Span}\, \{ \chi \, u_\infty \})$, provided that $p$
is chosen larger than (but close enough to) $\frac{N+1}{N-1}$. Let
us denote by $v_p$ this fixed point.

\medskip

Since $\delta > 1-N$, we have $|v_p| << u_{p}$ near $0$ and hence
the solution $u : = (1-\chi) \, \bar u_{p} + v_p$ is singular and
positive near $0$. We now prove that $u : = (1-\chi) \, \bar u_{p} +
v_p$ is also positive at infinity.  Indeed, the function $v_p$ can
be written as
\[
v_p = \tilde v_p + a_p \, \chi\, u_\infty
\]
where, according to formula \equ{formulaa}, $a_p$ can be computed as
\[
a_p \, (N-1) \, \int_{S^N_+} \,\theta_N\, d\sigma = \int_{{\mathbb
R}^N_+} (\Delta (1-\chi) \, \bar u_{p} + | (1-\chi) \, \bar u_{p} +
v_p|^{p} ) \, dx = \int_{{\mathbb R}^N_+} | (1-\chi) \, \bar u_{p} +
v_p|^{p} \, dx .
>0
\]
This implies that $a_p >0$, and by the maximum principle, it is now
easy to check that $u >0$ in ${\mathbb R}^N_+ -\{0\}$. This
completes the proof of Proposition~\ref{prop1}. \qed

\bigskip
\subsection{Some open questions}

The results  of proof of Propositions \ref{prop1} and \ref{prop2}
and their parallel with the radial case and $p$ close to $\frac
N{N-2}$, to which we recall ODE phase plane analysis applies, lead
us naturally to several questions concerning existence of solutions
of $\Delta u + u^p =0$ on the punctured half space ${\mathbb
R}^{N}_+-\{0\}$ with $0$ boundary data.    We list some of them
next.

\medskip

\noindent {\bf Question 1.} We believe that the  solution $ u_1$
which has been obtained in Proposition~\ref{prop1} for $p$ close to
${N+1\over N-1}$ should actually  exist for all $p \in
(\frac{N+1}{N-1}, \frac{N+2}{N-2})$.

\medskip
\noindent {\bf Question 2.} When $p = \frac{N+2}{N-2}$, we believe
that there exists a one parameter family of solutions of the form
\[
u (x) = |x|^{\frac{2-N}{2}} \, v ( - \log |x|,\theta)
\]
where $t \longmapsto v(t,  \cdot )$ is periodic. This one-parameter
family of solution corresponds to the well know periodic solutions
for the singular Yamabe problem and also to Delaunay surfaces in the
context of constant mean curvature surfaces.

\medskip

\noindent {\bf Question 3.} When  $p >\frac{N+2}{N-2}$, $N \geq 3$,
we believe that there exists a solution of $\Delta u + u^p=0$
defined on ${\mathbb R}^{N}_+$ which is identically equal to $0$ on
$\partial {\mathbb R}^N_+$ and which is asymptotic to $ u_0$ in
\equ{u00} at $\infty$. This solution should correspond to the smooth
radially symmetric solution of the same equation which is defined on
the whole space and decays like $|x|^{-\frac{2}{p-1}}$ at infinity,
when $p > \frac{N+2}{N-2}$.

\medskip

\noindent {\bf Question 4.}  Are there singular solutions when $p
\geq \frac{N+1}{N-3}$, $N \geq 4$? In this regime separation of
variables in general fails.

\medskip

Some partial answer to this question is given in \cite{BVPV}.

\bigskip
\setcounter{equation}{0}
\section{The bounded domain case:  proofs of Theorems \ref{th:1} and \ref{th:2}}

The proof of our main results relies on two basic ingredients: one
is the, already established, existence of the ``basic  cells'' given
by Propositions \ref{prop1} and \ref{prop2} which we will use to
construct approximations to singular solutions. Another important
ingredient, on which we elaborate in the next two subsections,  is
the analysis of invertibility of Laplace's operator, for right hand
sides exhibiting  a controlled singular behavior on a given embedded
submanifold of $\partial \Omega$, in the same spirit to that of
Lemma \ref{pr:3}. Then we will use a fixed point scheme analogous to
that in the proof of Proposition \ref{prop1}.

\medskip
For notational convenience, we will assume in what follows that
$\Omega$ is actually a  subset of $\RR^{N}$, and that $S$ is a
smooth embedded submanifold of $\partial \Omega\subset \RR^{N}$ with
dimension $k$. We define
$$
N= n+k. $$ We start by setting up a suitable description of the
space and Laplacian operator in natural coordinates associated to
$S$. While the analysis below is done for $k\ge 1$, it applies
equally well to the point-singularity case $k=0$, being actually
simpler.

\medskip
\subsection{Local coordinate system}
In a neighborhood of a point $\bar p$ of $S$ let us choose
coordinates $y_1, \ldots, y_k$ on $S$. Next we choose sections $E_1,
\ldots, E_{n-1}$ of the normal bundle of $S$ in $\partial \Omega$.
We can define Fermi coordinates in some tubular neighborhood of $S$
in $\partial \Omega$ by using the exponential map,
\[
F(p ,(x_1, \ldots, x_{n-1}))  =  \mbox{Exp}^{\partial \Omega}_p (
\sum_i x_i\, E_i(p))
\]
for $p$ in a neighborhood of $\bar p \in S$ and $(x_1, \ldots,
x_{n-1})$ in some neighborhood of $0$ in ${\mathbb R}^{n-1}$.

\medskip
In these coordinates, it is well known that the induced metric $g_b$
on $\partial \Omega$ can be expanded as
\[
g_b = g_{{\mathbb R}^{n-1}} + g_{S} + {\mathcal O} (|x|)
\]
where $g_S$ denotes the induced metric on $S$ and $x  = (x_1,
\ldots, x_{n-1})$.

\medskip
Finally, to parameterize a neighborhood of a point of $\partial
\Omega$ in $\Omega$, we denote by $E_n$ the normal (inward pointing)
vector field about $\partial \Omega$ and again use the exponential
map to define
\[
G(q ,x_n)  = q +  x_n \, E_n(q)
\]
for $q$ in a neighborhood of $\bar p$ in $\partial \Omega$ and $x_n
\geq 0$ in some neighborhood of $0$.

\medskip
In these coordinates, it is well known that the Euclidan metric in
$\Omega$ can be expanded as
\[
g_{{\mathbb R}^{n+k}} = dx_n^{2} + g_b + {\mathcal O}(x_n).
\]

\medskip
Collecting these two expansions, we conclude that in these
coordinates the (Euclidean) Laplacian  can be expanded as
\begin{equation}
\Delta  =  \Delta_{NS}  + {\mathcal O}(|x|) \, \nabla^{2} +
{\mathcal O}(1) \, \nabla \label{eq:3.14}
\end{equation}
where $x  = (x_1, \ldots, x_n)$ and $\Delta_{NS} =
\Delta_{{\mathbb R}^n} + \Delta_{S}$ is the Laplace Betrami
operator on $NS$, the normal bundle of $S$ in ${\mathbb R}^{n}$.

\subsection{Analysis of the Laplacian in weighted spaces}
We want to prove a result in the same spirit as that of Lemma
\ref{pr:3} in the current setting. To do this, we need to define
weighted spaces on $\bar \Omega \setminus S$, which have a
controlled blow up rate as $S$ is approached. Unlike those in Lemma
\ref{pr:3}, we choose H\"older spaces, which are more suitable to
deal with linear perturbations which are second order operators. Let
us define, for sufficiently small $R>0$, half ``balls'' and
``annuli''
\[
\bar B_+ (R): = \{ (p,x) \in NS_+ \  / \  |x| \in (0, R]\}
\]
and
\[
\bar A_+ (R_1, R_2) :=  \{ (p,x) \in NS_+ \ /\  |x|\in [R_1, R_2]\}
\]
$B_+ (R)$ is roughly ``half'' of a tubular neighborhood of radius
$R$ of the manifold $S$, or just a ball in case that $S$ reduces to
a single points. We consider the following weighted space of
functions defined on $\bar B_+ (R) \setminus S $.

\begin{definition}
The space ${\mathcal C}^{\ell, \alpha}_\delta (\bar B_+ (R)
\setminus S)$ is the space of functions $u \in {\mathcal C}^{\ell,
\alpha}_{loc} (\bar B_+ (R) \setminus S) )$ for which the norm
\[
\| u \|_{{\mathcal C}^{\ell, \alpha}_\delta( \bar B_+(R)\setminus S
)} = \sup_{r \in (0,R)} r^{-\delta} \, \| u( \cdot, r \, \cdot)
\|_{{\mathcal C}^{\ell, \alpha}( \bar A_+(r/2, r) ) }
\]
is finite.
\end{definition}

We consider now the problem
\begin{equation}
\begin{array}{rllllll}
  \Delta_{NS} u &=& |x|^{-2} f \quad \hbox{in $\bar B_+ (R) \setminus S $}
\\ [3mm]
  u &=& 0 \quad \hbox{on $\partial \bar B_+ (R) \setminus S $}.
\end{array}
\label{probl1}\end{equation} We have the validity of the following
result.

\begin{lemma} Assume that $\delta \in (1-n, 1)$.  There exists a constant
$c>0$, independent of $R>0$, such that, for each $ f\in {\mathcal
C}^{0, \alpha}_{\delta } (\bar B_+(R) \setminus  S)$, there is a
solution  $ u =  G_{\delta , R} \, (f) $ of problem $\equ{probl1}$
which defines a linear operator of $f$ and satisfies the estimate
$$
\| u \|_{{\mathcal C}^{2, \alpha}_\delta( \bar B_+(R)\setminus S )}
\, \le\, c \| f \|_{{\mathcal C}^{0, \alpha}_\delta( \bar
B_+(R)\setminus S )} \, .
$$
\label{pr:1}
\end{lemma}

\noindent{\bf Proof.} We only carry out the proof for $R=1$ since
the general case follows  by scaling. First we solve for each $r \in
(0, 1/2)$ the problem
\begin{equation}
\begin{array}{rllllll}
  \Delta_{NS} u &=& |x|^{-2} f \quad \hbox{in $A_+(r,1)                 $}
\\ [3mm]
  u &=& 0 \quad \hbox{on $\partial A_+(r,1)                 $}
\end{array}
\label{probl2}\end{equation} and call $u_r$ its unique solution.
Maximum principle employed in a similar way as in Lemma \ref{pr:3},
taking into account expansion (\ref{eq:3.14}), yield the {\it a
priori} bound
\[
|u_r| \leq c \, \| f \|_{{\mathcal C}^{0, \alpha}_\delta ( \bar B_+
(1) \setminus S )} \,  |x|^{\delta}
\]
where $c =c( n, \delta) >0$. Then, elliptic estimates applied on
geodesic balls of radius $r$ centered at distance $2r$ from $S$ give
the following  bound on the gradient of $u$
\[
|\nabla u_r |\leq c \, \| f \|_{{\mathcal C}^{0, \alpha}_\delta(
\bar B_+(1) \setminus S)} \, |x|^{\delta -1}
\]
for some $c = c(n, \delta) >0$. Using Arzela's theorem, we conclude
that, for a sequence of radii tending to $0$, the sequence $u_r$
converges to a function $u$ which satisfies
\[
|u|\leq c \,  \| f \|_{{\mathcal C}^{0, \alpha}_\delta( \bar B_+(1)
\setminus S )} \, |x|^{\delta}
\]
and solves \equ{probl1} for $R=1$. Again, elliptic estimates applied
on geodesic balls of radius $r$ centered at distance $2r$ from $S$
yield the bound
\[
\|  u \|_{{\mathcal C}^{2, \alpha}_{\delta } (\bar B_+ (1)
\setminus S)} \leq c \, \| f\|_{{\mathcal C}^{2, \alpha}_{\delta}
(\bar B_+ (1) \setminus S )}
\]
for some constant $c = c(n, \delta) >0$. Uniqueness of the limit $u$
is easy to get and we leave  it to the reader. The proof is
concluded. \hfill $\Box$

\bigskip
Next we will extend the previous result to the entire domain $\bar
\Omega \setminus  S$. To do so, we consider a function \[ \gamma :
\bar \Omega \setminus  S \longrightarrow (0, \infty)
\]
smooth, positive, which in the above defined local coordinates
coincides with $|x|$ in a neighborhood of $S$ in $\bar \Omega$. This
function will play the role of the function $|x|$ defined in $B_+(R)
\setminus  S$. We define accordingly weighted H\"older spaces as
follows.

\begin{definition}
We let the space ${\mathcal C}^{\ell, \alpha}_\delta (\bar \Omega
\setminus S )$ be that of functions $u \in {\mathcal C}^{\ell,
\alpha}_{loc} (\bar \Omega \setminus  S )$ for which the norm
\[
\| u\|_{{\mathcal C}^{\ell, \alpha}_\delta( \bar \Omega \setminus S
)} = \| u \|_{{\mathcal C}^{\ell, \alpha}_\delta (\bar B_+(R)
\setminus S )} + \| u \|_{{\mathcal C}^{\ell, \alpha} (\Omega
\setminus B_+(R/2))}
\]
is finite.
\end{definition}

We consider now the problem
\begin{equation}
\begin{array}{rllllll}
  \Delta u &=& \gamma^{-2} f \quad \hbox{in $\Omega \setminus S $}
\\ [3mm]
  u &=& 0 \quad \hbox{on $\partial \Omega \setminus S $} .
\end{array}
\label{probl3}\end{equation} We have the following result, extension
of Lemma \ref{pr:1}.

\begin{lemma} Assume that $\delta \in (1-n, 1)$.  There exists a constant
$c>0$ such that, for each $ f\in {\mathcal C}^{0, \alpha}_{\delta }
(\bar \Omega \setminus  S)$, there is a solution  $ u = G_{\delta }
\, (f) $ of problem  $\equ{probl3}$ which defines a linear operator
of $f$ and satisfies the estimate
$$
\| u \|_{{\mathcal C}^{2, \alpha}_\delta( \bar \Omega \setminus S )}
\, \le\, c \| f \|_{{\mathcal C}^{0, \alpha}_\delta( \bar \Omega
\setminus S )} \, .
$$
\label{pr:2}
\end{lemma}

\noindent{\bf Proof.} The proof follows from Lemma~\ref{pr:1},
expansion (\ref{eq:3.14}) and a linear perturbation argument. First,
we claim that the result of Lemma~\ref{pr:1} remains true in $\bar
B_+(R) \setminus  S$ if the operator $\Delta_{NS}$ is replaced by
$\Delta$ and if $R$ is chosen small enough. Indeed, we have from
(\ref{eq:3.14}) and Proposition~\ref{pr:1}
\[
\| f -  \gamma^{2} \, (\Delta - \Delta_{NS}) \circ G_{\delta, R} (f)
\|_{{\mathcal C}^{0, \alpha}_{\delta} (\bar B_+ (R) \setminus S)} \leq c \,
R \, \| f\|_{{\mathcal C}^{0, \alpha}_{\delta} (\bar B_+ (R) \setminus S)} .
\]
The claim follows at once from a perturbation argument, provided
that $R$ is fixed small enough. We denote by $\bar G_{\delta , R}$
the right inverse for $\Delta$ in $\bar B_+(R) \setminus  S$.

\medskip

We consider a cut-off function $\chi_R$ which is equal to $1$ in
$B_+(R/2) \setminus S$ and equal to $0$ in $\bar \Omega
\setminus B_+(R)$. We define
\[
\tilde f  : = f - \gamma^{2} \, \Delta (\chi_R u_1),
\]
where $u_1 =   \bar G_{\delta,R} (f)$. Observe that this function
is supported in $\bar \Omega \setminus B_+(R/2)$. We have that
$\tilde f \in {\mathcal C}^{0, \alpha}(\bar \Omega)$ and
\[
\| \tilde f  \|_{{\mathcal C}^{0, \alpha} (\bar \Omega)} \leq c \|
f\|_{{\mathcal C}^{0, \alpha}_{\delta} (\bar \Omega \setminus
\Gamma )}
\]
for some constant $c = c(n, \delta, R) >0$.

\medskip

Finally, we can solve
\[
\begin{array}{rllllll}
  \Delta u_2 &=& \gamma^{-2} \tilde f \quad \hbox{in $\Omega $}
\\ [3mm]
  u_2 &=& 0 \quad \hbox{on $\partial \Omega  $} .
\end{array}
\]
We have the bound
\[ \| u_2 \|_{{\mathcal C}^{2, \alpha}  (\bar \Omega  )} \leq c \|
\tilde f\|_{{\mathcal C}^{0, \alpha}(\bar \Omega \setminus S )} .
\] The desired result then follows by letting the solution of
\equ{probl3} be $u = u_1 +u_2$. \hfill $\Box$

\bigskip
\subsection{Proof of Theorems~\ref{th:1} and \ref{th:2}}

We are now in a position to provide the proof of Theorem~\ref{th:1}
and Theorem~\ref{th:2}. The argument goes along the same lines as
that in the proof of Proposition \ref{prop1}, now with Lemma
\ref{pr:2} playing the role of Lemma \ref{pr:3}.

\medskip
 We recall that  we are now
assuming that $\Omega$ is a domain in $\RR^{N}$. We also write
\[
N = n+k .\]

\medskip
\noindent\noindent{\bf Proof of Theorem \ref{th:1}, case $p=\frac{n
+1}{n-1}$.} We assume that  $S$ is either a finite number of points
of $\partial \Omega$, namely $k=0$, or an embedded $k$-dimensional
submanifold of $\partial \Omega$. For all $\varepsilon >0$ small
enough, we define
\[
u_\varepsilon  : = \chi_R  \, \varepsilon^{n-1} \, u_1(\varepsilon
\,|x|, \ve x_n)
\]
where $u_1$ is the solution provided by Proposition~\ref{th:2} and
$\chi_R$ is a cut-off function which equals $1$ in $B_+ (R)
\setminus S$ and  $0$ in $\Omega \setminus B_+ (2R)$. Here we fix $R
>0$ sufficiently small and use for $x$ and $x_n$ the meanings
given in the previous subsections. In particular we have that $u_\ve
=0$ on $\partial\Omega \setminus S$.

\medskip
The problem we want to solve then reads
\begin{equation}
\begin{array}{rlllll}
\Delta(u_\varepsilon + v) + | u_{\varepsilon} +v|^{\frac{n+1}{n-1}}
&=&0 \quad \hbox{in $\Omega  $,}\\ [3mm]  v &=&0 \quad \hbox{on
$\partial\Omega\setminus S $,} \nonumber\end{array}
\nonumber\end{equation} where we also require $u_\varepsilon + v
>0$ in $\Omega$. Let us fix $\delta \in (1-n, 2-n]$. By virtue of
Lemma~\ref{pr:2}, we can rewrite this equation as the fixed point
problem
\begin{equation}
v  = - G_\delta \, \left( \gamma^2 (\Delta u_\varepsilon + |
u_{\varepsilon} +v|^{\frac{n+1}{n-1}} ) \right) .
\label{fp4}\end{equation}

We have the validity of the following fact: there is a constant $c_0
= c(\delta , \Omega , S)
>0$ such that
\[
\|\gamma^{2} ( \Delta u_\varepsilon +
u_\varepsilon^{\frac{n+1}{n-1}} )\|_{{\mathcal C}^{0,
\alpha}_{\delta } (\bar \Omega \setminus S)} \leq c_0 \, (\log
(1/\varepsilon) )^{\frac{1-n}{2}},
\]
result that is a consequence of  expansion (\ref{eq:3.14}) and a
direct computation using the asymptotic properties of $u_1$ in
Proposition \ref{prop2}.

\medskip
We  restrict our attention to the case where $\delta \leq 2-n$ since
$\gamma^2 \, (\Delta u_\varepsilon +
u_\varepsilon^{\frac{n+1}{n-1}})$ is  bounded by a constant times
$|x|^{2-n}$ near $S$, and $\delta \leq 2-n$ guarantees that this
function belongs to ${\mathcal C}^{0, \alpha}_{\delta} (\bar \Omega
\setminus S)$.

\bigskip
A second estimate we can directly check is  the following: Assume
that $\delta  \in (1-n,  2-n]$ is fixed. There exists a constant $c
= c(\delta , \Omega ,S) >0$ such that
\[
\|\gamma^{2} ( |u_\varepsilon + v_2|^{\frac{n+1}{n-1}} -
|u_\varepsilon + v_1|^{\frac{n+1}{n-1}} )\|_{{\mathcal C}^{0,
\alpha}_{\delta} (\bar \Omega \setminus S)} \leq c \, (\log (1/\varepsilon)
)^{-1} \, \| v_2 -v_1\|_{{\mathcal C}^{2, \alpha}_{\delta} (\bar
\Omega \setminus S)}
\]
for all $v_2, v_1 \in {\mathcal C}^{2, \alpha}_{\delta} (\bar \Omega
\setminus S)$ satisfying
\[
\| v_i\|_{{\mathcal C}^{2, \alpha}_{\delta} (\bar \Omega \setminus S)} \leq
\, 2 \, c_0 \, (\log (1/\varepsilon) )^{\frac{1-n}{2}} .
\]

The above estimates allow an application of contraction mapping
principle in the ball of radius $2 \, c_0 \, (\log (1/\varepsilon)
)^{\frac{1-n}{2}}$ in ${\mathcal C}^{2, \alpha}_{\delta} (\bar
\Omega \setminus S)$ to predict existence of a solution to problem
\equ{fp4}, which we denote by $v_\varepsilon$.

\medskip
Since $\delta > 1-n$, we have $|v_\varepsilon| << u_\varepsilon$
near $S$ and hence the solution $u : = u_\varepsilon +
v_\varepsilon$ is singular along $S$ and is positive near $S$. The
maximum principle then implies that $u >0$ in $\Omega$. This
completes the proof of Theorem~\ref{th:1} in the case
$p=\frac{n+1}{n-1}$. \qed

\bigskip
{\bf The proof of  Theorem~\ref{th:2}, case $p=\frac{n+1}{n-1}$.}
This proof  uses similar arguments together with an induction
process. By assumption, $A$ is closed and contains a sequence of
$k$-dimensional submanifolds $S_i$, $i\in {\mathbb N}$ such that
$\cup_i S_i$ is dense in $A$. We define inductively the sequence of
functions $u_i$ which are solutions of
\begin{equation}
\Delta u_i + u_i^{\frac{n+1}{n-1}} =0 \label{eq:6.6}
\end{equation}
in $\Omega$, satisfy $u_i=0$ on $\partial \Omega - \cup_{j=0}^i S_j$
and are singular along $\cup_{j=0}^i  S_j$. Assume for example that
$u_{i-1}$ has already been constructed, then, we define
\[
\tilde u_{i} = u_{i-1} + \varepsilon^{n-1}_{i} \, \chi_{r_{i+1}} \,
u_1 (\varepsilon_{i} \, \mbox{dist}(\cdot, S_i))
\]
where $u_1$ is the solution provided by Proposition~\ref{th:2},
$r_i$ is fixed small enough less than half the distance from $S_i$
to $\cup_{j=0}^{i-1} S_j$ and $\varepsilon_i
>0$ is small enough. Applying a perturbation argument as above, we
can perturb $\tilde u_i$ into a solution $u_i = \tilde u_i + v_i$ of
(\ref{eq:6.6}) for some function $v_i \in {\mathcal C}^{2,
\alpha}_{\delta} (\bar \Omega \setminus \cup_{j=0}^i S_j)$. Taking
$\varepsilon_i$ small enough, we can ensure that
\begin{equation}
\| u_i  - u_{i-1} \|_{L^1 (\Omega)} \leq 2^{-i} \label{eq:AA}
\end{equation}
\begin{equation}
\| \mbox{dist}(\cdot, \partial \Omega)^2 \, | u_i  -
u_{i-1}|^{\frac{n+1}{n-1}} \|_{L^1 (\Omega)}^{\frac{n+1}{n-1}} \leq
2^{-i} \label{eq:BB}
\end{equation}
and
\begin{equation}
\| \tilde \gamma^{\delta} v_i \|_{L^\infty(\Omega)} \leq 2^{-i}
\label{eq:CC}
\end{equation}
where $\tilde \gamma = \mbox{dist} (\cdot, S)$ and where $\delta \in
(1-n, 2-n]$ is fixed. Clearly (\ref{eq:AA}) ensures that the
sequence $(u_i)_i$ converges in $L^1(\Omega)$ to a function $u$.
Moreover (\ref{eq:AA}) and (\ref{eq:BB}) imply that $u$ is a weak
solution of (\ref{eq:1}). Finally, (\ref{eq:CC}) implies that the
nontangential limit of $u$ at any point of $S$ is equal to
$+\infty$. \qed

\medskip
Finally, we observe that, using Proposition~\ref{prop1} instead of
Proposition~\ref{prop2}, the results of Theorems ~\ref{th:1} and
\ref{th:2} hold when $p
> \frac{n+1}{n-1}$, is sufficiently close to $\frac{n+1}{n-1}$. The only
difference being that, in the proof of the result corresponding to
the one of Theorem~\ref{th:2}, in addition to the properties
(\ref{eq:AA}) to (\ref{eq:CC}) which ensure the convergence of the
sequence of solutions in the appropriate space, we may also ask that
the sequence converges in $W^{1, q}(\Omega)$, for some $q$ close
enough to $1$.  The proofs are concluded. \qed

\bigskip
\centerline{\bf Acknowledgement } This work has been supported by
grants Ecos/Conicyt C05E05, Fondecyt 1030840, 104936, and  FONDAP.

\noindent \textsc{M. del Pino}

\medskip

\noindent Departamento de Ingenier\'{\i}a Matem\'atica and CMM,
Universidad de Chile, Casilla 170 Correo 3, Santiago, Chile.

\noindent {\bf email :} delpino@dim.uchile.cl,

\medskip

\noindent \textsc{M. Musso}

\medskip

\noindent Departamento de Matem\'atica, Pontificia Universidad
Cat\'olica de Chile, Avda. Vicu\~na Mackenna 4860, Macul, Chile and
Dipartimento di Matematica, Politecnico di Torino, Corso Duca degli
Abruzzi 24, 10129 Torino, Italy.

\noindent {\bf email :} mmusso@mat.puc.cl,

\medskip

\noindent \textsc{F. Pacard}

\medskip

\noindent Universit\'e Paris 12 and Institut Universitaire de France

\noindent  {\bf email :} pacard@univ-paris12.fr
\end{document}